%% file: root.tex
\documentclass[
nolists,
nofoot,
subscriptcorrection,
upint,
barcolor=RedOrange!80,
mathalfa=cal=cm,
balance,
hyphenate,
]{asmejour} %

\input{preamble.tex}

\begin{document}

\SetAuthorBlock{Ruairi Moran\,$^{\ast,\,2}$\\
Sheila Bagley\,$^{\ast\ast}$\\
Seth Kasmann\,$^{\ast\ast}$\\
Rob Martin\,$^{\ast\ast}$\\
David Pasley\,$^{\ast\ast}$\\
Shane Trimble\,$^{\ast\ast}$\\
James Dianics\,$^{\ast\ast}$\\
\vspace{0.2em}
Pantelis Sopasakis\,$^{\ast}$}{
\vspace{1em}
$^{\ast}$\,Queen's University Belfast,\\
EEECS, i-AMS Centre,\\
Ashby Building, BT9 5AH, Belfast, UK\\
\vspace{1em}
$^{\ast\ast}$\,EquipmentShare,\\
5710 Bull Run Dr, Columbia, MO, 65201, USA
} 

\title{NMPC for Collision Avoidance by Superellipsoid Separation$^{1}$}

\keywords{Collision avoidance, Superellipsoids, Nonlinear MPC}

\begin{abstract}
\input{abstract.tex}
\end{abstract}

\date{
$^{1}$\,Paper presented at the 2024 Modeling, Estimation, and Control Conference (MECC 2024), Chicago, IL, Oct. 28-30, 2024. Paper No. MECC2024-63.\\
\hspace*{4mm}$^{2}$\,Corresponding author, email: \texttt{rmoran05@qub.ac.uk}
}%% This command must come before \maketitle

\maketitle %% Essential!

\section{Introduction}\label{sec:introduction}
In the rapidly evolving domains
of robotics and autonomous systems,
two critical motion challenges stand out:
path planning and trajectory
planning~\citep{khan2021comprehensive}.
Path planning involves finding a global path through complex
environments that a mobile robot can follow.
Trajectory planning involves finding a local path between
points of the global path, while
avoiding unexpected obstacles.
This paper explores trajectory planning.
Achieving safe and efficient
obstacle avoidance is a challenge
that spans various engineering applications,
from autonomous ground vehicles
(AGVs)~\citep{liu2018nonlinear}
and
drones~\citep{xue2021vision}
to industrial
robots~\citep{gai20196}
and many more~\citep{patle2019review}.
This work is applied to
skid-steer heavy equipment
in a construction site
environment~\citep{iii1995autonomous, melenbrink2020site},
where ensuring operator safety is
paramount~\citep{teizer2010autonomous}.
Skid-steer vehicles
have many applications due to their
simple mechanical structure
and high mobility~\citep{khan2021comprehensive}.

To address the challenge of collision-free 
control, we explore
the optimization-based approach of a
nonlinear model predictive controller
(NMPC)~\citep{grune2017nonlinear}.
An NMPC calculates safe
and feasible trajectories
and control inputs for a vehicle in real time,
through a combination
of nonlinear model-based prediction
and constrained optimization.
This paper proposes a
collision avoidance condition for an NMPC,
derived from the combination of
the separating hyperplane theorem (SHT)
\cite[Sec. 2.5.1]{boyd2004convex}
and superellipsoids.

Superellipsoids, a generalization
of ellipsoids and rectangles,
offer versatile object representation.
Superellipsoids have been used
to represent objects for obstacle
avoidance, 
however, to only represent the obstacles~\citep{pauls2022real, menon2017trajectory},
to only represent the vehicle~\citep{smith2017implementing},
or, for path planning using a semi-infinite constrained optimization problem rather than optimal control~\citep{wang2004real}.

In this paper we implement a navigation system as
a two-layer
hierarchical NMPC~\citep{scattolini2009architectures},
with a high-level NMPC layer
that involves the
proposed collision avoidance condition,
and a low-level NMPC layer
that tracks the reference trajectory from
the high-level layer.
Hierarchical NMPC
enables different sampling times
and control horizons
for each layer,
ensuring quick and adaptive responses
to environmental changes
at the low layer
while allowing optimal
trajectory planning at the high layer.

\begin{figure}
    \centering
    \includegraphics[width=\linewidth]{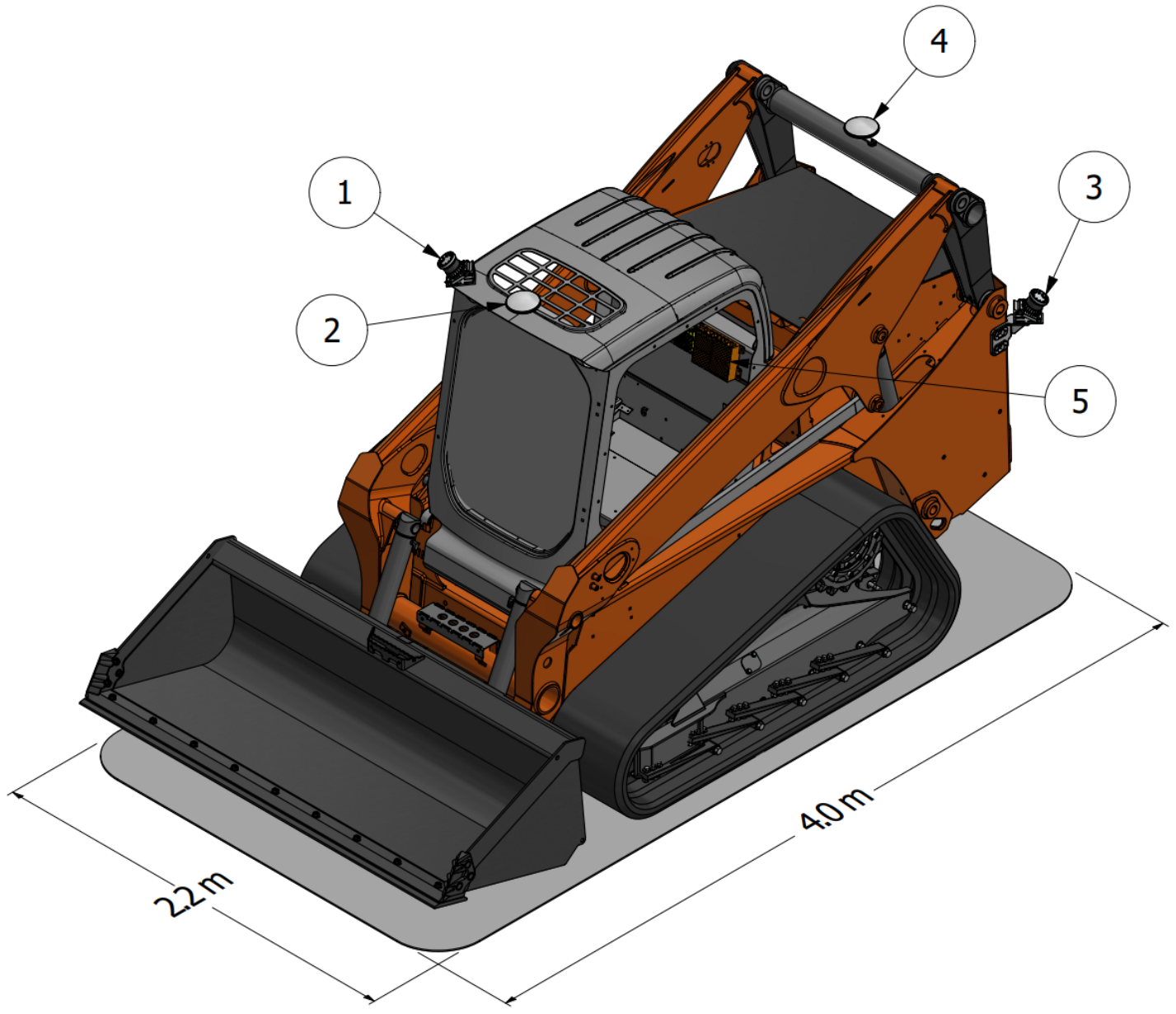}
    \caption{EquipmentShare's
        compact track loader,
        \emph{Moe}.
        The labels are:
        1. front LiDAR,
        2. front GNSS,
        3. rear LiDAR,
        4. rear GNSS,
        5. IMU and TX2.}
    \label{fig:moe}
\end{figure}

Note that generally mobile robot
optimal control
techniques for collision avoidance model the
vehicle as a point
or an ellipsoid,
and obstacles as ellipsoids or
polygons~\citep{khan2022design, rafai2022review}.
Such approaches can naturally lead to overly
conservative overapproximation, thus necessitating
a large vehicle-to-obstacle distance.
We contribute a less conservative
approach that models both vehicle and obstacles
as superellipsoids in an NMPC problem formulation that allows
its fast solution using OpEn,
and demonstrate the trajectory-planning capability with a real-time application.

\subsection{Notation}
Let $\N_{[t_1,t_2]}$ denote the set of natural numbers from $t_{1}$ to $t_{2}$,
$\R$ the set of real numbers,
and $\R^{n}$ the $n$-dimensional real space equipped with
the standard inner product $\langle \cdot, \cdot \rangle$.
Then,
the support function
of a set $\mc{X}\subseteq\R^{n}$ at $x\in\R^{n}$ is defined by
${
    \supp_{\mc{X}}(x) = \Sup_{y\in \mc{X}} \langle x, y \rangle,
}$
and the $p$-norm of $x$
is defined by
${
    \|x\|_{p}{}={}\left(
    \sum_{i=1}^{n} |x_{i}|^{p}
    \right)^{\frac{1}{p}},
}$
where $p\in\left[ 1, \infty \right)$.
We denote the
unit $p$-norm ball
\begin{equation}
    \B_{p} {}={}
    \left\{
    x \in\R^{n} :
    \|x\|_{p} \leq 1
    \right\}.
\end{equation}
Lastly, we define the
matrix-set product
${
    A\mc{X} {}={}
    \left\{
    Ax, \forall x\in\mc{X}
    \right\}
}$.

\section{Vehicle}\label{sec:system-dynamics}
The vehicle in this paper is
the compact track loader (CTL) \emph{Moe},
illustrated in Fig.~\ref{fig:moe},
which is a type of skid-steer that moves
on two parallel tracks.
The arm and bucket are held in
place in this work.

The global frame of reference
is the North-East-Down (NED) coordinate
system~\cite[Sec. 2.2]{cai2011unmanned}.
We denote a position on
the North axis by $x_{1}$
and on the East axis by $x_{2}$.
The body-fixed frame of reference
uses the same coordinate system
fixed to the vehicle's principal axes.
The vehicle's heading is denoted $\theta$,
where $\theta = 0$ corresponds to the
vehicle facing North.
A positive heading
corresponds to a clockwise
rotation of the vehicle
with respect to the global frame of reference.

A CTL has nonholonomic dynamics and 
steers by rotating the tracks at different rates.
The vehicle's state $z=(c, \theta, v)$
consists of
the position $c\in\R^{2}$,
heading $\theta\in (-\pi,\pi]$,
and orbital velocity $v\in\R$.
The vehicle's input $u=(r, s)$ consists of
the throttle $r \in [-1, 1]$
and spin $s \in [-1, 1]$.
We use the unicycle model,
which is a simplified differential
drive~\cite[Sec. 13.2.4.1]{lavalle2006planning}
\begin{subequations}\label{eq:dynamics-ct}
    \begin{align}
        \dot{c}_{1}
         & {}={}
        v \cos(\theta),
        \\
        \dot{c}_{2}
         & {}={}
        v \sin(\theta),
        \\
        \dot{\theta}
         & {}={}
        \alpha s,
        \\
        \dot{v}
         & {}={}
        \beta
        \left(
        r v_{\rm max} - v
        \right),
    \end{align}
\end{subequations}
where $\alpha, \beta, v_{\rm max} \geq 0$ are tuning parameters
for adapting the model to
different terrains
(softer ground will result in
smaller $\alpha$ and $\beta$)
and different vehicles
($v_{\rm max}$ is the maximum velocity).
Despite its simplicity, the model is sufficient for
trajectory planning and control tasks
as demonstrated in Ref.~\cite[Sec. 4.4]{khan2022design}
especially provided the vehicle motion is slow 
as shown in Ref.~\citep{pazderski2008trajectory}.

The CTL uses a diesel-hydraulic powertrain to drive each track separately.
The sensors on board include
an iNEMO ISM330DLC
IMU updating at $208\, \unit{\hertz}$,
two Swift Navigation Piksi Multi GNSS
units at $10\, \unit{\hertz}$,
and two
LiDAR units at $20\, \unit{\hertz}$.
These provide readings
with standard measurement error of
the vehicle's latitude and longitude
to $10\, \unit{\milli\meter}$,
altitude
to $15\, \unit{\milli\meter}$,
heading
to $0.6\unit{\degree}$,
and obstacle detection within a $35\, \unit{\meter}$ radius.

The sensors are connected via local Ethernet to an NVIDIA Jetson TX2,
which is the main computational unit on the vehicle.
The TX2 is connected to
the vehicle control unit (VCU) over
controller area network (CAN) bus,
which allows the TX2
to override
the original equipment manufacturer's (OEM)
throttle and spin inputs.
The on-board equipment is depicted
in Fig.~\ref{fig:moe}.

The TX2 codebase is primarily
written in Python.
The vehicle state and control inputs
are logged at $10\, \unit{\hertz}$.
The controllers in
Section~\ref{sec:solver-implementation}
are built with Optimization Engine
(OpEn)~\citep{sopasakis2020open},
an open-source code generation software
for embedded nonlinear optimization
that runs in Rust.
OpEn was chosen because
it has low memory requirements,
involves simple algorithm operations,
the generated Rust code is provably memory safe,
and has been shown to be faster than
interior point and
sequential quadratic programming
implementations~\citep{agv2018}.

\section{Collision Avoidance Formulation}\label{sec:constraint-formulation}
While this work focuses on collision avoidance for a CTL, this section can be applied to general collision avoidance problems.

\subsection{Introducing Superellipsoids}\label{sec:intro-superellipsoids}
We will be modeling the vehicle and obstacles
by sets of the form
\begin{equation}\label{eq:superellipsoid-mink-sum}
    \mc{X} {}={}
    \left\{
    RSx + c \,\vert\,
    x\in\B_{p}
    \right\},
\end{equation}
where
$p \in [2, \infty)$,
$x\in\R^{2}$,
$c \in\R^{2}$ is the center point,
$S \in\R^{2 \times 2}$ is a diagonal scaling matrix,
$S = \operatorname{diag}(s_1, s_2)$,
where $s_{1},s_{2}>0$,
and $R_{\theta}$, 
that belongs to the orthogonal group ${\rm SO}(2)$,
is the elementary rotation matrix
\begin{equation}
    R_{\theta} =
    \begin{bmatrix}
        \cos\theta & -\sin\theta \\
        \sin\theta & \cos\theta
    \end{bmatrix},
\end{equation}
where $\theta\in (-\pi,\pi]$.
We restrict $p \geq 2$
where the sets are smooth and convex.
When $p=2$, Eq.~\eqref{eq:superellipsoid-mink-sum}
describes an ellipsoid.
The matrix product $RS$
scales the set $\B_{p}$
by $s_{1}$ on the North axis,
by $s_{2}$ on the East axis,
and rotates the scaled set clockwise
about the center point
through the angle $\theta$.
Note that
$R$ is orthogonal
and
$RS$ is nonsingular.

Superellipsoid is
the umbrella term
for a set described by
Eq.~\eqref{eq:superellipsoid-mink-sum}~\cite[p.1760]{weisstein2002crc}.
The set is also called
a Lam\'e curve~\citep{gridgeman1970lame}. 
In this paper we use superellipsoids with $p=3$
(see Fig.~\ref{fig:sep-hyp}). Higher values of $p$
will make the shape look closer to a rectangle.

\subsection{Collision Detection}
Let us denote the vehicle set $\mc{V}$
and an obstacle set $\mc{E}$.
Hereafter, we will only consider
one obstacle, however, the extension
to multiple obstacles will be obvious.

A collision between two sets
can be detected by
a condition that checks
if the sets are disjoint.
We know from the separating hyperplane
theorem~\cite[Sec. 2.5.1]{boyd2004convex}
that
if two convex sets $\mc{V}$ and $\mc{E}$
are disjoint,
there exists
a separating hyperplane
such that $\mc{V}$ is on one side
and $\mc{E}$ is on the other.
There are several converse SHTs,
such as the separating axis theorem
(SAT)~\cite[Sec. 5.2.1]{ericson2004real}.
We propose Theorem~\ref{th:strict-separation-3}.
\begin{theorem}[Converse SHT]
    \label{th:strict-separation-3}
    Let $\mc{V}, \mc{E} \subseteq \R^{n}$
    be nonempty sets, there exists
    $a\in\R^{n}$,
    and
    \begin{equation}\label{eq:theorem}
        \Sup_{x \in\mc{V}}\ \langle a, x \rangle
        {}<{}
        \Inf_{x \in\mc{E}}\ \langle a, x \rangle.
    \end{equation}
    Then $\mc{V}$ and $\mc{E}$
    are disjoint,
    i.e., $\mc{V} \cap \mc{E}=\emptyset$.
\end{theorem}
\textit{Proof.}
Suppose Eq.~\eqref{eq:theorem} holds,
but $\mc{V}$ and $\mc{E}$
are not disjoint, i.e.,
$\mc{V} \cap \mc{E}\neq\emptyset$.
Take $y \in \mc{V} \cap \mc{E}$.
Then
\(
\Sup_{x \in\mc{V}}\ \langle a, x \rangle
{}\geq{}
\langle a, y \rangle
{}\geq{}
\Inf_{x \in\mc{E}}\ \langle a, x \rangle,
\)
which contradicts Eq.~\eqref{eq:theorem}. $\Box$

\begin{figure}
    \centering
    \includegraphics[width=\linewidth]{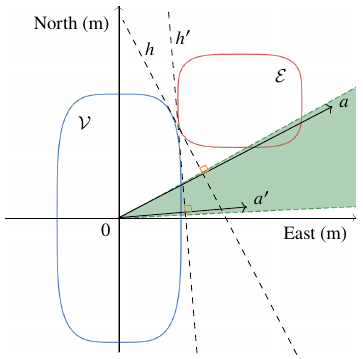}
    \caption{\mathversion{normal}
        Examples of separating axes 
        are labeled $a$ and $a'$. 
        Disjoint sets $\mc{V}$ and $\mc{E}$
        lie in opposite open halfspaces
        on either side of the corresponding
        dashed black hyperplanes, 
        labeled $h$ and $h'$ respectively.
        The green open set
        shows all the separating axes
        that satisfy the condition
        in Eq.~\eqref{eq:theorem}
        for this scenario.
        Note that the green set is 
        closed under positive scalar multiplication.}
    \label{fig:sep-hyp}
\end{figure}

In this paper,
an axis $a$ that satisfies Eq.~\eqref{eq:theorem}
is called a separating axis.
Figure~\ref{fig:sep-hyp} illustrates
two edge-case examples of separating axes for the same scenario.
The condition
can be equivalently written as
\begin{equation}\label{eq:sep-hyp-the}
    \supp_{\mc{V}}(a)
    + \supp_{-\mc{E}}(a)
    {}<{}
    0.
\end{equation}
Lastly, note that for safety reasons, the vehicle superellipsoid
is enlarged, therefore, the strict inequality in Eq.~\eqref{eq:sep-hyp-the} can be relaxed to the inequality $\leq$.

\subsection{Collision Avoidance Condition}
The vehicle and obstacles are modeled by
closed superellipsoids.
Let us define the vehicle
\(
\mc{V} {}={}
\left\{
R^{v}_{\theta} S^{v} x + c^{v} \,\vert\,
x\in\B_{p}
\right\},
\)
and an obstacle as
\(
\mc{E} {}={}
\left\{
R^{e}_{\theta} S^{e} x + c^{e} \,\vert\,
x\in\B_{p}
\right\}.
\)
The support function of $\mc{V}$ at $a$ is
\begin{align}
    \supp_{\mc{V}}(a)
     & {}={}
    \Sup_{\substack{w \in R^{v}_{\theta} S^{v} \B_{p} + c^{v}}}
    \langle a, w \rangle
    \notag
    \\
     & {}={}
    \Sup_{d \in R^{v}_{\theta} S^{v} \B_{p}}
    \langle a, d + c^{v} \rangle
    \text{, where } d = w - c^{v}
    \notag
    \\
     & {}={}
    \supp_{R^{v}_{\theta} S^{v} \B_{p}}(a)
    + \langle a, c^{v} \rangle.
\end{align}
Similarily, as the sets
$R_{\theta}S\B_{p}$ and $-R_{\theta}S\B_{p}$
are equal, the support function
of $-\mc{E}$ at $a$ is
\begin{align}
    \supp_{-\mc{E}}(a)
    {}={}
    \supp_{R^{e}_{\theta} S^{e} \B_{p}}(a)
    + \langle a, -c^{e} \rangle.
\end{align}
The support function of a superellipsoid
centered at the origin
is the dual of the $p$-norm~\cite[Sec. 7.2]{bauschke2017convex}
\(
\supp_{R_{\theta}S\B_{p}}(a)
{}={}
\|S R_{\theta}^\tr a\|_{q},
\)
where $p\in\left[ 2, \infty \right)$,
and $p$ and $q$ are conjugate
exponents,
i.e.,
$\frac{1}{p} + \frac{1}{q} = 1$.
The pair $(p,q)=(1,\infty)$
is valid and
$
    \|x\|_{\infty}
    =
    \max_{i}
    |x_{i}|
$.
Now, we can express Eq.~\eqref{eq:sep-hyp-the} as
\begin{equation}\label{eq:strict-sep-set-constraint-no-epsilon}
    \|S^{v} R^{v\tr}_{\theta} a\|_{q}
    + \|S^{e} R^{e\tr}_{\theta} a\|_{q}
    + \langle a, c^{v} - c^{e} \rangle
    < 0.
\end{equation}
The set of all $a$ that satisfy
the condition in
Eq.~\eqref{eq:strict-sep-set-constraint-no-epsilon}
is closed under positive scalar multiplication,
therefore, we can assume that ${\|a\|_{2}=1}$.

\section{Hierarchical NMPC Design}\label{sec:solver-implementation}
A two-layer
hierarchical NMPC is implemented to control
the vehicle.
While this work focuses on a CTL, this two-layer controller can be applied to any slow-moving dual-tracked vehicle.
The high-level
trajectory-planning
controller is denoted HC and
the low-level
trajectory-following
controller is denoted LC.
Collision avoidance is only enforced in HC, as the heavy equipment is assumed to be moving slowly such that a collision constraint violation would not be possible in between HC time steps. In addition, LC is assumed to closely follow the trajectory output of HC. Until HC produces a new trajectory, LC will continue to follow the last trajectory provided by HC. Lastly, we assume static obstacles described by superellipsoids.

The input to HC includes
the vehicle state,
the target state,
and any obstacles.
The output is a list
of the vehicle's state
at each predicted stage.
The input to LC includes
the trajectory output of HC.
The output of LC is the
throttle and spin inputs
for the vehicle
to follow the trajectory.
In line with the receding horizon control
strategy, the first control action
from LC is applied to the vehicle's inputs.

To formulate the NMPC problems,
we discretize Eq.~\eqref{eq:dynamics-ct}
by forward Euler discretization,
$f(z_{t},u_{t}) = z_{t} + T \dot{z}_{t}$,
with sampling time $T > 0$.
Then, we define $F_{0}(z,u) = z$
and recursively
$
    F_{i+1}(z,u)
    = f\left(
    F_{i}(z,u), u
    \right)
$.

Next, we present the hierarchical
NMPC, where the sampling time of LC
is $t_{L}=T$ and 
the sampling time of HC is $\tau$ times longer,
$t_{H}=\tau T$,
to allow enough time
for the high-level optimization problem
to be solved.

\subsection{High-Level Problem}\label{sec:high-level-problem}
We formulate the
high-level problem with
the vehicle state
$z_{t}=(c_{t}, \theta_{t}, v_{t})$
and vehicle input $u_{t}=(r_{t}, s_{t})$
at stage $t$,
as defined in Section~\ref{sec:system-dynamics},
and the system dynamics
$F_{\tau}$.
The high-level discrete-time
trajectory-planning
optimal control problem is
\begin{subequations}\label{eq:solver-problem-high}
    \begin{align}
        \mathbb{P}_{H}: \hspace*{-0.8em} & \minimise_{\substack{
                (u_t)_{t=0}^{H-1},\
        (z_t)_{t=0}^{H},                                                  \\
                (a_t)_{t=0}^{H}}}
        \sum_{t=0}^{H-1}\ell_{t}(z_t, u_t) + V_f(z_{H}),
        \\
        \textbf{s.t.}\
                                         & z_{t+1} = F_{\tau}(z_t, u_t),\
        t\in\N_{[0, H-1]},
        \\
                                         &
        |r_{t}| \leq r_{\rm max},\
        t\in\N_{[0, H-1]},\label{eq:ph-throttle}
        \\
                                         &
        |s_{t}| \leq s_{\rm max},\
        t\in\N_{[0, H-1]},\label{eq:ph-spin}
        \\
                                         &
        \label{eq:separating-constraint}
        \|(S^{v} R^{v\tr}_{\theta})_{t} a_{t}\|_{q}
        + \|(S^{e} R^{e\tr}_{\theta})_{t} a_{t}\|_{q}
        \notag
        \\
                                         & \qquad
        + \langle a_{t}, c^{v}_{t} - c^{e}_{t} \rangle
        \leq
        0,
        t\in\N_{[0, H]},
        \\
                                         &
        \|a_{t}\|_{2} = 1,\
        t\in\N_{[0, H]},
        \label{eq:axis-constraint}
        \\
                                         &
        z_0 = z,
    \end{align}
\end{subequations}
for some
$p \geq 2$,
and $1 \leq \tau\in\N$,
where
$q=\frac{p}{p-1}$,
$\ell_{t}$ is the stage cost function at stage $t$,
$V_{f}$ is the terminal cost function,
and $H\in\N$ is the horizon.

We can extend $\mathbb{P}_{H}$ to consider
multiple obstacles by adding multiple constraints of the forms of~\eqref{eq:separating-constraint} and \eqref{eq:axis-constraint}---one of each constraint for each obstacle.
The size of $\mathbb{P}_{H}$ increases with the number of obstacles and is limited by the computational resources available.

To reduce the problem size,
$\ell_{t}\left(z_{t}, u_{t}\right)=0$
for all odd $t$,
and for all even $t$
the stage cost function is
\begin{align*}
     & \ell_{t}\left(z_{t}, u_{t}\right)
    {}={}
    q_{c} \|c_{t} - c_{\rm ref}\|^{2}_{2}
    +
    q_{\theta} \left(\theta_{t} - \theta_{\rm ref}\right)^{2}
    \notag
    \\
     & \quad
    +
    q_{r} r_{t}^{2}
    +
    q_{r_{\Delta}} \left(r_{t} - r_{t-1}\right)^{2}
    +
    q_{s} s_{t}^{2}
    +
    q_{s_{\Delta}} \left(s_{t} - s_{t-1}\right)^{2},
\end{align*}
and the terminal cost function is
\begin{equation*}
    V_{f}\left(z_{H}\right)
    {}={}
    q_{c_{H}} \|c_{H} - c_{\rm ref}\|_{2}^{2}
    +
    q_{\theta_{H}} \left(\theta_{H} - \theta_{\rm ref}\right)^{2},
\end{equation*}
where
$
    q_{c},
    q_{\theta},
    q_{r},
    q_{r_{\Delta}},
    q_{s},
    q_{s_{\Delta}},
    q_{c_{H}},
    q_{\theta_{H}}
    > 0$ are the weights,
$c_{\rm ref}$ is the reference position
and $\theta_{\rm ref}$ is the reference heading.
The weight of $q_{c}$ and $q_{c_{H}}$ penalize
the Euclidean distance between the vehicle and the target,
$q_{\theta}$ and $q_{\theta_{H}}$ penalize the heading difference
between the vehicle and the target,
$q_{r}$ and $q_{s}$ penalize the use of control inputs,
and $q_{r_{\Delta}}$ and $q_{s_{\Delta}}$ penalize changes in the control inputs.

\subsection{Low-Level Problem}
The solution $(z^{\star}_t)_{t=0}^{H}$
of $\mathbb{P}_{H}$,
where
$z_{t}^{\star}=(c_{t}^{\star}, \theta_{t}^{\star}, v_{t}^{\star})$,
is a trajectory
towards the target state,
comprising
states every $t_{H}$ seconds
for LC to follow.
The trajectory is tracked
by an NMPC that penalizes
the difference between the
predicted state at stage $k$
and the next state
in the planned trajectory
after $k$, at stage
$t_{k}=\max\left\{1, \left\lceil \tau^{-1} k \right\rceil\right\}$,
where
$\lceil \cdot \rceil$
is the ceiling function.

The trajectory tracking is
similar to the
two-and-a-half carrots
method in Ref.~\citep{reiter2015two},
where the most important predicted state
is at the beginning
rather than the end.
In LC, the stage of the most important state
is denoted $\omega$ and is selected
so that
$\omega$ is a multiple of $\tau$,
is close to
the beginning of the planned trajectory,
and the trajectory is updated before
the vehicle has time to reach it,
e.g., $\omega = 2\tau$.
This results in good trajectory-following performance
as the vehicle prioritizes
following the trajectory immediately before it,
rather than the vehicle's state
at the end of the trajectory.

We differentiate notation that appear
in both problems by an overline
on the low-level variables and functions.
We formulate the low-level problem with
the vehicle state
$\bar{z}_{k}=(\bar{c}_{k}, \bar{\theta}_{k}, \bar{v}_{k})$
and vehicle input
$\bar{u}_{k}=(\bar{r}_{k}, \bar{u}_{k})$
at stage $k$,
as defined in Section~\ref{sec:system-dynamics},
and the system dynamics
$f$.
The low-level discrete-time
trajectory-following
optimal control problem is
\begin{subequations}\label{eq:solver-problem-low}
    \begin{align}
        \mathbb{P}_{L}:
        \minimise_{\substack{
        (\bar{u}_k)_{k=0}^{L-1}, \\
                (\bar{z}_k)_{k=0}^{L}}}
         &
        \sum_{\substack{k=0      \\k\neq\omega}}^{L-1}
        \bar{\ell}(\bar{z}_{k}, \bar{u}_k)
        + V_{\omega}(\bar{z}_{\omega},\bar{u}_{\omega})
        + \bar{V}_f(\bar{z}_{L}),
        \\
        \textbf{s.t.}\
         &
        \bar{z}_{k+1} = f(\bar{z}_k, \bar{u}_k),\
        k\in\N_{[0, L-1]},
        \\
         &
        |\bar{r}_{k}| \leq r_{\rm max},\
        k\in\N_{[0, L-1]},\label{eq:pl-throttle}
        \\
         &
        |\bar{s}_{k}| \leq s_{\rm max},\
        k\in\N_{[0, L-1]},\label{eq:pl-spin}
        \\
         &
        \bar{z}_0 = \bar{z},
    \end{align}
\end{subequations}
for some $\omega\in\N_{[0,L-1]}$,
where
$\bar{\ell}$ is the stage cost function,
$V_{\omega}$ is a high-weighted stage cost function,
$\bar{V}_{f}$ is the terminal cost function,
and $L\in\N$ is the horizon.
The horizon time of
$\mathbb{P}_{H}$
must be longer than that of
$\mathbb{P}_{L}$
in order for the planned trajectory
to have enough states
to follow,
i.e., $Lt_{L} < Ht_{H}$.

The stage cost function is
\begin{align}
     & \bar{\ell}\left(\bar{z}_{k}, \bar{u}_{k}\right)
    {}={}
    \bar{q}_{c} \|\bar{c}_{k} - c_{t_{k}}^{\star}\|^{2}_{2}
    +
    \bar{q}_{\theta} \left(\bar{\theta}_{k} - \theta_{t_{k}}^{\star}\right)^{2}
    \notag
    \\
     & \quad
    +
    \bar{q}_{r} \bar{r}_{k}^{2}
    +
    \bar{q}_{r_{\Delta}} \left(\bar{r}_{k} - \bar{r}_{k-1}\right)^{2}
    \notag
    +
    \bar{q}_{s} \bar{s}_{k}^{2}
    +
    \bar{q}_{s_{\Delta}} \left(\bar{s}_{k} - \bar{s}_{k-1}\right)^{2},
\end{align}
the high-weighted stage cost function is
\begin{align}
     & V_{\omega}\left(\bar{z}_{\omega}, \bar{u}_{\omega}\right)
    {}={}
    \bar{q}_{c_{\omega}} \|\bar{c}_{\omega} - c_{t_{k}}^{\star}\|_{2}^{2}
    +
    \bar{q}_{\theta_{\omega}} \left(\bar{\theta}_{\omega} - \theta_{t_{k}}^{\star}\right)^{2}
    \notag
    \\
     & \quad
    +
    \bar{q}_{r} \bar{r}_{k}^{2}
    +
    \bar{q}_{r_{\Delta}} \left(\bar{r}_{k} - \bar{r}_{k-1}\right)^{2}
    \notag
    +
    \bar{q}_{s} \bar{s}_{k}^{2}
    +
    \bar{q}_{s_{\Delta}} \left(\bar{s}_{k} - \bar{s}_{k-1}\right)^{2},
\end{align}
and the terminal cost function is
\begin{equation}
    \bar{V}_{f}\left(\bar{z}_{L}\right)
    {}={}
    \bar{q}_{c_{L}} \|\bar{c}_{\omega} - c_{t_{k}}^{\star}\|_{2}^{2}
    +
    \bar{q}_{\theta_{L}} \left(\bar{\theta}_{\omega} - \theta_{t_{k}}^{\star}\right)^{2},
\end{equation}
where
$
    \bar{q}_{c},
    \bar{q}_{\theta},
    \bar{q}_{r},
    \bar{q}_{r_{\Delta}},
    \bar{q}_{s},
    \bar{q}_{s_{\Delta}},
    \bar{q}_{c_{\omega}},
    \bar{q}_{\theta_{\omega}}
    \bar{q}_{c_{L}},
    \bar{q}_{\theta_{L}}
    > 0$ are the weights.
The weight of $\bar{q}_{c},\bar{q}_{c_{\omega}}$ and $\bar{q}_{c_{L}}$ penalize
the Euclidean distance between the vehicle and the target,
$\bar{q}_{\theta},\bar{q}_{\theta_{\omega}}$ and $\bar{q}_{\theta_{L}}$ penalize the heading difference,
$\bar{q}_{r}$ and $\bar{q}_{s}$ penalize the use of control inputs,
and $\bar{q}_{r_{\Delta}}$ and $\bar{q}_{s_{\Delta}}$ penalize changes in the control inputs.

\subsection{Tuning}\label{sec:nmpc-tuning}
To avoid local minima when solving $\mathbb{P}_H$,
we have two solvers:
one solver is enabled constantly and
we warm-start this solver's separating axes
with the repeated vector of the vehicle's current center
to the obstacle's center, that is,
${
    (c^{e} - c^{v}_{0})_{t},
    t\in\N_{[0, H]},
}$
where $c^{v}_{0}$ is the vehicle's center at stage $t=0$.
This is inspired by the approach in
Ref.~\cite[Sec. 5.2.1]{ericson2004real}
for polyhedral objects.
The other solver is enabled when
a previous solution that converged is available.
When both solvers are enabled,
they run in parallel
and the output solution of
the solver that converges with
the lowest cost
is used.
Then, the solution trajectory is accepted
if it has an adequately low
measure of infeasibility$^{3}$\footnotetext[3]{
    That is, if the infinity-norm of succesive estimates
    of the Lagrange multipliers is less than
    a tolerance ($0.001$) times the current penalty.
    This is a measure of the distance between
    $G_{1}(\xi, \bm{u})$ and $\mathcal{C}$,
    i.e., a measure of infeasibility of
    Eq.~\eqref{eq:alm}.
},
which is provided by the OpEn solver~\citep{sopasakis2020open}.
If not, the vehicle continues on the
last accepted trajectory.

The parameters for HC and LC are recorded in Table~\ref{tab:parameters}. We used experimental data
to estimate the demonstration model parameters from
Section~\ref{sec:system-dynamics}.

\begin{table}[t]
\caption{Simulation and demonstration parameters for HC and LC}
\begin{center}
\label{tab:parameters}
\begin{tabular}{p{12mm} l l|p{12mm} l l}
\toprule
HC & & & LC & & \\
\midrule
Param. & Sim. & Demo. & Param. & Sim. & Demo. \\
\midrule
$\alpha$ & 1 & 0.4 & $\alpha$ & 1 & 0.4 \\
$\beta$ & 0.2 & 0.15 & $\beta$ & 0.2 & 0.15 \\
$v_{\rm{max}}\ (\unit{\meter/\second})$ & 1 & 1 & $v_{\rm{max}}\ (\unit{\meter/\second})$ & 1 & 1 \\
$r_{\rm{max}}$ & 1 & 0.5 & $r_{\rm{max}}$ & 1 & 0.5 \\
$s_{\rm{max}}$ & 1 & 0.2 & $s_{\rm{max}}$ & 1 & 0.2 \\
$q_{c}$ & 1 & 100 & $\bar{q}_{c}$ & 100 & 100 \\
$q_{\theta}$ & 0 & 0 & $\bar{q}_{\theta}$ & 0 & 0 \\
$q_{r}$ & 0.01 & 0.01 & $\bar{q}_{r}$ & 0.01 & 0.1 \\
$q_{s}$ & 0.5 & 7 & $\bar{q}_{s}$ & 0.1 & 0.1 \\
$q_{r_{\Delta}}$ & 0 & 0 & $\bar{q}_{r_{\Delta}}$ & 0 & 5 \\
$q_{s_{\Delta}}$ & 0 & 0.01 & $\bar{q}_{s_{\Delta}}$ & 0 & 5 \\
& & & $\bar{q}_{c_{\omega}}$ & 1000 & 1000 \\
& & & $\bar{q}_{\theta_{\omega}}$ & 0 & 0 \\
$q_{c_{H}}$ & 20 & 5000 & $\bar{q}_{c_{L}}$ & 100 & 100 \\
$q_{\theta_{H}}$ & 0 & 0 & $\bar{q}_{\theta_{L}}$ & 0 & 0 \\
$t_{H}$ & 1 & 1 & $t_{L}$ & 0.1 & 0.1\\
$H$ & 40 & 40 & $L$ & 100 & 100 \\
& & & $T$ & 0.1 & 0.1 \\
& & & $\tau$ & 10 & 10 \\
& & & $\omega$ & 20 & 20 \\
\bottomrule
\end{tabular}
\end{center}
\end{table}

\subsection{Embedded Numerical Optimization}\label{sec:building}
OpEn
solves parametric nonconvex problems
of the form
\begin{subequations}\label{eq:open}
    \begin{align}
        \hspace{-1em}
        \mathbb{P}
        \left(
        \xi
        \right):
        \minimise_{\bm{u}\in\R^{n}}\
         & g(\xi, \bm{u}),
        \\
        \textbf{s.t.}\
         & \bm{u}\in\mathcal{U},
        \label{eq:proj}
        \\
         & G_{1}(\xi, \bm{u}) \in \mathcal{C},
        \label{eq:alm}
        \\
         & G_{2}(\xi, \bm{u}) = 0,
        \label{eq:pm}
    \end{align}
\end{subequations}
where $\bm{u}$
is the decision variable
and $\xi$
is a parameter vector.
Briefly, the requirements are that
$g$ is a smooth function with Lipschitz gradient,
$G_{1}$ and $G_{2}$ satisfy certain regularity conditions,
$\mathcal{U}$ is closed
and can be projected onto,
and the set $\mathcal{C}$ is
closed, convex, and
the point-to-set distance can be
computed~\citep{sopasakis2020open}.

Within the algorithm of OpEn,
constraints of the type in Eq.~\eqref{eq:proj}
are imposed by
projecting onto $\mc{U}$,
\eqref{eq:alm} by using
the augmented Lagrangian method (ALM),
and \eqref{eq:pm} by using
the quadratic penalty method (PM).

In particular,
to solve $\mathbb{P}_{H}$
we first eliminate the sequence of states
using the single shooting approach~\citep{sopasakis2020open},
that carries out the minimization
over the sequence of control actions.
Then, we use Eq.~\eqref{eq:proj}
to cast Constraints~\eqref{eq:ph-throttle},
\eqref{eq:ph-spin},
and \eqref{eq:axis-constraint},
and \eqref{eq:alm}
for \eqref{eq:separating-constraint}.
Note that while Constraint~\eqref{eq:separating-constraint}
can be cast in the form of
both Eq.~\eqref{eq:alm}
and \eqref{eq:pm},
the ALM significantly outperforms
the PM in this case.

Regarding $\mathbb{P}_{L}$,
we use the single shooting formulation again
and \eqref{eq:pl-throttle}
and \eqref{eq:pl-spin} are
projectable constraints of the
form in Eq.~\eqref{eq:proj}.

\section{Simulation}
\begin{figure}
    \centering
    \includegraphics[width=\linewidth]{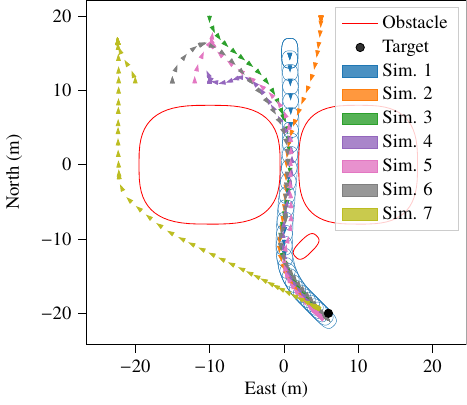}
    \caption{Seven simulated trajectories of a vehicle. The vehicle's and obstacle's parameters are described in Tables~\ref{tab:sim-obstacles-params} and \ref{tab:sim-vehicle-params}. The blue superellipsoids illustrate the vehicle size, which is the same for all simulations.}
    \label{fig:sim-path-plot}
\end{figure}

Hereafter, all measurements are given in meters unless otherwise stated, and all positions are (N, E).
In Fig.~\ref{fig:sim-path-plot} we present seven simulations that were run on
an NVIDIA Jetson TX2 fitted for \textit{Moe}.
All objects were modeled with $p=3$.
For all simulations, 
the vehicle was modeled
by a superellipsoid with parameters
$s_{1}=2.0$ and $s_{2}=1.1$,
the target point was $(-20, 6)$,
and the obstacle's states are described in Table~\ref{tab:sim-obstacles-params}.
For each simulation, the vehicle's initial states are described in Table~\ref{tab:sim-vehicle-params}.

\begin{table}[t]
\caption{Simulation parameters of the three obstacles}
\begin{center}
\label{tab:sim-obstacles-params}
\begin{tabular}{l c c c c c}
\toprule
Obstacle & $c_{1}$ & $c_{2}$ & $\theta\ (\unit{\radian})$ & $s_{1}$ & $s_{2}$ \\
\midrule
East & 0 & 10 & 0 & 8 & 8 \\
West & 0 & -10 & 0 & 8 & 9.5 \\
South & -11 & 3 & -0.79 & 2 & 1 \\
\bottomrule
\end{tabular}
\end{center}
\end{table}

\begin{table}[t]
\caption{Simulation parameters of the vehicle's initial states for seven simulations}
\begin{center}
\label{tab:sim-vehicle-params}
\begin{tabular}{l c c c c c c c}
\toprule
Param. & 1 & 2 & 3 & 4 & 5 & 6 & 7 \\
\midrule
$c_{1}$ & 15 & 20 & 20 & 11 & 11 & 11 & 11 \\
$c_{2}$ & 0.8 & 5 & -10 & -10 & -12 & -15 & -20 \\
$\theta\ (\unit{\radian})$ & 3.14 & 3.14 & 3.14 & 0 & 0 & 0 & 0 \\
\bottomrule
\end{tabular}
\end{center}
\end{table}

Figure~\ref{fig:sim-box-plot} illustrates
the runtimes of the high- and low-level solvers
during simulations. The low-level controller
had a median execution time of $1.1\, \unit{\milli\second}$
(max $13.9\, \unit{\milli\second}$), and the high-level
controller had a median of $434.5\, \unit{\milli\second}$
(max $714.0\, \unit{\milli\second}$).

\begin{figure}
    \centering
    \input{seven_sim_box_plot.tex}
    \caption{High- and low-level solver runtime data collected during simulations (see Fig.~\ref{fig:sim-path-plot})}
    \label{fig:sim-box-plot}
\end{figure}

\section{Demonstration}
A demonstration of the hierarchical controller
on one of EquipmentShare's CTLs was run
at their R\&D facility.
The vehicle autonomously navigates
around a virtual obstacle
to the target.
The demonstration is available online
at \texttt{https://youtu.be/HLNVZtbsl0c}.
The logged vehicle states
are plotted in Fig.~\ref{fig:demo-plot}
and the logged vehicle control inputs
are plotted in Fig.~\ref{fig:demo-inputs}.
The tracking error between the vehicle position
and the path planned by HC
is plotted in
Fig.~\ref{fig:demo-path-follow},
which shows LC closely
follows the trajectory planned by HC.

Again, all positions are (N, E).
The initial vehicle position was (-103.68, 77.51)
and the target position was (-79.50, 61.40).
The virtual obstacle was modeled
by a superellipsoid with parameters
$c_{1}=-92.2, c_{2}=69.8,
    s_{1}=1.0, s_{2}=1.5,$
and
$\theta=-0.5\, \unit{\radian}$.
The vehicle was modeled
by a superellipsoid with parameters
$s_{1}=2.0$ and $s_{2}=1.1$.
Both were modeled with $p=3$.
\begin{figure}
    \centering
    \input{demo_plot.tex}
    \caption{Vehicle state
        data collected
        during demonstration.
        The demonstration ends
        when the vehicle reaches
        within a meter of the target position.
        The blue superellipsoids depict the
        vehicle,
        growing darker as time evolves.
    }
    \label{fig:demo-plot}
\end{figure}
\begin{figure}
    \centering
    \input{i_inputs_plot.tex}
    \caption{Throttle and spin
        control actions
        from LC
        during the demonstration.
        The vehicle initially
        considers turning left
        to go around the obstacle,
        then decides on a trajectory
        to the right of the obstacle.%
    }%
    \label{fig:demo-inputs}
\end{figure}

\begin{figure}
    \centering
    \input{tracking_error.tex}
    \caption{Tracking error between the vehicle position and the planned path provided by HC.
    This shows that
    the vehicle position closely
    follows the planned trajectory,
    while the trajectory is updated. The error had a 95th percentile distance of $42\, \rm{mm}$ (max $80\, \rm{mm}$).
    Recall that the GNSS unit has a standard meaurement error of $10\, \rm{mm}$ for these readings.
    The max compound error of $90\, \rm{mm}$
    is the main safety reason
    for enlarging the superellipsoids.}
    \label{fig:demo-path-follow}
\end{figure}

During the demonstration,
the solver in HC
that converged with the lowest cost
had a 95th percentile runtime of
$185.7\, \unit{\milli\second}$
with a maximum of
$900\, \unit{\milli\second}$,
and the solver in
LC had a 95th percentile runtime of
$6.4\, \unit{\milli\second}$
with a maximum of
$90\, \unit{\milli\second}$.
The maximum runtimes are
90\% of the solver's sampling time.
As a safety measure,
if LC reaches its maximum runtime,
the control input is set to zero
so the vehicle does not move.
If HC reaches its maximum runtime,
LC continues to track
the last trajectory provided by HC.
Neither solver reached its maximum runtime
during the demonstration.

\section{Conclusions and Future Work}\label{sec:conclusion}
We developed
a novel approach to
obstacle avoidance in NMPC,
leveraging a combination of
the separating hyperplane theorem
and superellipsoids.
We proposed a tractable
collision avoidance condition
and demonstrated the efficacy
and real-time capability of the approach
through simulations and experimental results.

The use of superellipsoids for
object representation
has shown to be advantageous
in comparison to modeling
with ellipsoids,
as rectangular geometries
that are often encountered
in real-world scenarios
---
such as vehicles and buildings
---
can be more accurately described
with the same complexity.

Future work will focus on extending the approach to handle dynamic and nonconvex obstacles in more challenging environments, and on a more thorough theoretical analysis to provide stability and feasibility guarantees for the high- and low-level solvers.

\section*{Funding Data}
\begin{itemize}
\item This paper has received funding by the project ``BotDozer: GPU-accelerated
model predictive control for autonomous heavy equipment,'' which is part of the Doctoral Training Grant No S3809ASA, funded by EPSRC.
\end{itemize}

\bibliographystyle{asmejour}
\bibliography{database.bib}

\end{document}

%% file: preamble.tex
\usepackage{bm}
\usepackage{amsmath,amsfonts,mathtools}
\usepackage[noend]{algorithmic}
\usepackage{graphicx}
\usepackage{nicefrac}
\usepackage{siunitx}
\usepackage{romannum}
\usepackage{enumitem}

\fancypagestyle{title}{%
\fancyhf{} % clear all header and footer fields
}
\fancypagestyle{plain}{%
\fancyhf{} % clear all header and footer fields
}

\usepackage{pgfplots}
\usepackage{tikz}
\usepgfplotslibrary{fillbetween}
\usetikzlibrary{decorations.pathmorphing,
    patterns,shapes,arrows,positioning, backgrounds, intersections, calc}
\tikzstyle{Nodes}=[circle, draw=black, fill=orange!20, line width=1.5pt, minimum size=15pt]
\tikzstyle{arrow} = [line width=1pt,->,>=stealth]
\tikzstyle{axis} = [line width=1pt,->,>=stealth]
\usepackage{tkz-euclide}
\usetikzlibrary{external}
\usetikzlibrary{pgfplots.groupplots}
\tikzstyle{block}=[
draw=white,
thick,
minimum width=1cm,
minimum height=1cm,
inner sep=4pt,
text width = 7em,
text centered,
fill = white,
]

\hypersetup{
	pdfauthor={Ruairi Moran},                       		   	% <=== change to YOUR name[s]!
	pdftitle={ASME LDSC Paper},                  	% <=== change to YOUR pdf file title
	pdfkeywords={ASME LDSC paper, LaTeX, BibTeX style, asmejour class},% <=== change to YOUR pdf keywords
	pdfsubject = {ASME LDSC paper},			% <=== change to YOUR subject
}

\newcommand{\tr}{{\intercal}}
\newcommand{\R}{{\rm I\!R}}
\newcommand{\N}{{\rm I\!N}}

\newcommand{\B}{\mathcal{B}}

\newcommand{\supp}{\delta^{*}}

\newcommand{\mc}[1]{\mathcal{#1}}

\DeclareMathOperator*{\Sup}{\textbf{sup}}
\DeclareMathOperator*{\Inf}{\textbf{inf}}

\DeclareMathOperator*{\minimise}{\textbf{Minimize}}

\definecolor{myred}{rgb}{0.8,0.0,0.0}
\definecolor{mygreen}{rgb}{0.0,0.6,0.0}
\definecolor{myblue}{rgb}{0.0,0.0,0.8}

\newtheorem{theorem}{Theorem}

%% file: abstract.tex
This paper introduces a novel
NMPC formulation for real-time obstacle avoidance
on heavy equipment
by modeling both vehicle and obstacles
as convex superellipsoids.
The combination of this approach with
the separating hyperplane theorem
and Optimization Engine (OpEn)
allows to achieve efficient obstacle avoidance
in autonomous heavy equipment and robotics.
We demonstrate the efficacy of the approach
through simulated and experimental results,
showcasing a skid-steer loader's capability to navigate
in obstructed environments.

%% file: seven_sim_box_plot.tex
% This file was created with tikzplotlib v0.10.1.
\begin{tikzpicture}

\definecolor{darkgray176}{RGB}{176,176,176}
\definecolor{lightgray}{RGB}{211,211,211}

\begin{groupplot}[group style={group size=1 by 2}]
\nextgroupplot[
height=3cm,
tick align=outside,
tick pos=left,
width=8cm,
x grid style={lightgray},
xlabel near ticks,
xmajorgrids,
xmin=0, xmax=800,
xtick style={color=black},
y grid style={darkgray176},
ylabel near ticks,
ymin=-54, ymax=56,
ytick style={color=black},
ytick={1},
yticklabels={High}
]
\addplot [black]
table {%
353.75 -49
353.75 51
508 51
508 -49
353.75 -49
};
\addplot [black]
table {%
353.75 1
171 1
};
\addplot [black]
table {%
508 1
714 1
};
\addplot [black]
table {%
171 -24
171 26
};
\addplot [black]
table {%
714 -24
714 26
};
\addplot [black, opacity=0.8, mark=+, mark size=3, mark options={solid,fill opacity=0,draw=red}, only marks]
table {%
63 1
};
\addplot [red]
table {%
434.5 -49
434.5 51
};

\nextgroupplot[
height=3cm,
tick align=outside,
tick pos=left,
width=8cm,
x grid style={lightgray},
xlabel near ticks,
xlabel={Runtime (ms)},
xmajorgrids,
xmin=0, xmax=20,
xtick style={color=black},
y grid style={darkgray176},
ylabel near ticks,
ymin=-4.5, ymax=6.5,
ytick style={color=black},
ytick={1},
yticklabels={Low}
]
\addplot [black]
table {%
0.28275 -4
0.28275 6
4.2625 6
4.2625 -4
0.28275 -4
};
\addplot [black]
table {%
0.28275 1
-0.498 1
};
\addplot [black]
table {%
4.2625 1
9.894 1
};
\addplot [black]
table {%
-0.498 -1.5
-0.498 3.5
};
\addplot [black]
table {%
9.894 -1.5
9.894 3.5
};
\addplot [black, opacity=0.8, mark=+, mark size=3, mark options={solid,fill opacity=0,draw=red}, only marks]
table {%
13.469 1
13.938 1
};
\addplot [red]
table {%
1.1255 -4
1.1255 6
};
\end{groupplot}

\end{tikzpicture}

%% file: tracking_error.tex
% This file was created with tikzplotlib v0.10.1.
\begin{tikzpicture}

\definecolor{darkgray176}{RGB}{176,176,176}

\begin{axis}[
height=3cm,
tick align=outside,
tick pos=left,
width=8cm,
x grid style={darkgray176},
xlabel={Time (s)},
xmin=-0.337586629390717, xmax=106.357249724865,
xtick style={color=black},
y grid style={darkgray176},
ylabel={Error (mm)},
ymin=-3.8030206598699, ymax=83.6676067715659,
ytick style={color=black}
]
\addplot [semithick, black]
table {%
4.51217865943909 11.8512007963858
4.62650728225708 10.4205351535804
4.68955993652344 0.232085892108736
4.89270973205566 3.41914159382466
5.30114459991455 11.9388624860129
5.69852113723755 25.2176701873186
6.09428954124451 16.3634592067324
6.53304553031921 30.6058343940174
6.90561127662659 28.9015790113762
7.2950451374054 33.4328959659453
7.69099807739258 31.274482370508
8.08839726448059 31.7550633082922
8.4898738861084 33.4339359891449
8.88959121704102 25.9355482578205
9.30010414123535 28.9365287062241
9.68976640701294 30.3051411781029
10.1026804447174 25.7736991398739
10.4876523017883 35.404109813348
10.8927955627441 32.0190013505817
11.2855000495911 20.8239886082411
11.6858358383179 35.2786408471379
12.084135055542 31.9781626952287
12.4822182655334 34.8226322605063
12.8819973468781 38.9393688804653
13.2807958126068 26.4994715661502
13.6777913570404 2.29375672720915
14.0906069278717 6.01347658028715
14.4897172451019 8.98145321918924
14.8827471733093 13.4379336362846
15.28253865242 13.2854201527723
15.6785070896149 17.9577204703168
16.0650537014008 11.9826815654079
16.4634208679199 20.780478386778
16.8723692893982 28.5572476801995
17.2696702480316 27.9049551458397
17.6616427898407 27.7066774121643
18.0599050521851 26.7905277634789
18.4593210220337 29.6071710304766
18.878585100174 12.5491513602209
19.2571213245392 5.66239238505287
19.6735513210297 18.716211744074
20.0569679737091 8.81082797021732
20.4681923389435 22.7469039957055
20.862832069397 7.88158012988751
21.2727830410004 6.46513937976055
21.6707949638367 3.74492838257231
22.066300868988 3.46736904973532
22.4589977264404 10.0725329158861
22.8598003387451 3.41265747746926
23.2509973049164 6.42262070980901
23.6558840274811 3.13065762250258
24.0572454929352 11.5136179525691
24.4490375518799 27.0401960898765
24.8471217155457 18.6126651172788
25.2555527687073 18.5512205252072
25.6463572978973 10.9466045088978
26.0589334964752 1.92748692798331
26.4439375400543 1.56836316716022
26.8465433120728 4.04139630641767
27.2496821880341 1.7550140353532
27.6475048065186 0.456227962132608
28.0409705638885 12.4308734745327
28.4550323486328 19.0477882804542
28.8640582561493 8.28150923566048
29.2583417892456 4.46899796885964
29.6531887054443 3.11913907760367
30.0530230998993 3.85255219316679
30.4655458927155 1.63799942559268
30.8522877693176 14.958515763092
31.2534482479095 2.90251003997169
31.6519782543182 0.172916950649912
32.051561832428 8.55214910396798
32.4471364021301 17.69899872766
32.8586778640747 24.554359563915
33.2452518939972 24.7071701070629
33.6624464988708 33.7353173316378
34.051766872406 37.5158414243236
34.4442026615143 46.4540699725978
34.8663935661316 45.1319430995959
35.2764256000519 55.7144222135672
35.6472206115723 56.8639797472855
36.044921875 79.6916691610461
36.4519317150116 13.9996337617291
36.8448767662048 2.82506312867561
37.2413382530212 2.95656094564021
37.655818939209 0.841423467441703
38.0415456295013 7.95695348402835
38.4488949775696 2.60748224008831
38.8441686630249 8.51832490894519
39.2482447624207 18.3779858490183
39.6534848213196 24.4795178628275
40.0554196834564 13.2966287503382
40.4549922943115 16.4524314460695
40.8435795307159 7.36059849944084
41.2449057102203 0.600517021924586
41.6394855976105 10.1466318873864
42.0781154632568 6.02227012369006
42.4415187835693 10.1050515026276
42.8465387821198 15.6167215168955
43.2457115650177 34.9474029236791
43.6345882415771 46.9385775469529
44.0632936954498 39.5085693860337
44.4439907073975 17.2075817965632
44.8324906826019 9.33136429660147
45.2488408088684 15.512385805296
45.6301038265228 37.0936210183596
46.0302639007568 4.66149783691804
46.4537448883057 1.65146217798613
46.838850736618 29.036508765919
47.2334899902344 22.8613174745708
47.6236252784729 15.7676020315188
48.0325922966003 5.08332948149827
48.4236736297607 11.2493701663245
48.8182611465454 28.5850335793375
49.2307333946228 39.4450558224719
49.6230380535126 15.9886047965169
50.0230073928833 5.0044068068661
50.4326162338257 19.558714338996
50.8457028865814 39.183431853191
51.2481484413147 31.3003187905637
51.620566368103 25.7347096439952
52.0173764228821 17.8233049309993
52.4296875 10.5993174106729
52.8155839443207 34.3027238766368
53.2131283283234 42.2181440299143
53.615091085434 52.8423101680946
54.0141553878784 3.90130852075695
54.4289586544037 8.79797098444719
54.8121297359467 17.955235827823
55.2102696895599 30.8649894781534
55.6130774021149 50.2255173732151
56.0163986682892 11.8088338784401
56.4325680732727 9.15498812199592
56.8325099945068 16.9276705579422
57.2268435955048 5.04120508207093
57.6177217960358 9.48458251289961
58.021538734436 24.409920772612
58.4048516750336 1.14981169575433
58.8045749664307 5.00440287792351
59.2111370563507 3.20025898877224
59.6106402873993 3.81608089969672
59.9978926181793 16.3001953388194
60.4076900482178 1.08507055819207
60.8114535808563 6.3626082137858
61.2020773887634 25.5694136866155
61.6035482883453 31.253769877794
62.0003910064697 11.1392849906125
62.403623342514 5.88372233049411
62.7912828922272 7.6789749390336
63.1897881031036 21.3396263464588
63.6028709411621 28.7134385395236
63.9980351924896 12.7415304457682
64.3953502178192 15.4253779057604
64.7972157001495 13.470092716636
65.1834444999695 37.129545116041
65.5817184448242 33.0953718408308
65.9801487922668 19.1860432195792
66.3884217739105 24.5488169201737
66.7955701351166 26.3420203984615
67.1914343833923 32.5593096093029
67.5906972885132 43.0404710039754
67.988970041275 40.1109701937436
68.392951965332 1.24218726297481
68.7765355110168 5.09012385322693
69.1803736686707 13.6470541949037
69.5980813503265 30.8466593271207
69.9820845127106 44.5008844765369
70.382985830307 5.59707057633884
70.7777981758118 4.68534709730734
71.1700179576874 14.8015121444401
71.576827287674 20.3581945760559
71.9795739650726 31.1574563607343
72.368914604187 18.3103249306409
72.7797491550446 24.5425321709526
73.1660594940186 27.7060639764799
73.5657060146332 32.3878605979528
73.971337556839 39.9339858557766
74.3687374591827 10.8605179803896
74.7671558856964 19.0107048125789
75.1585576534271 20.2477333917147
75.5874006748199 21.9326475557916
75.9817881584167 25.6440306957777
76.3583557605743 29.0215055001387
76.7564392089844 0.667837958022824
77.1652269363403 7.69835860863381
77.5564067363739 9.163998068622
77.9698278903961 7.5486247609595
78.3786492347717 17.0813111216414
78.7537717819214 8.27009507356339
79.1636538505554 5.56198759532328
79.564866065979 13.4989493569877
79.9484713077545 11.0173779629636
80.359189748764 14.1653703481209
80.7429163455963 3.77485994406266
81.1581406593323 8.07131759311241
81.5573964118958 7.46486522726993
81.9446918964386 11.454686883747
82.3845138549805 7.42870030500642
82.7452375888824 6.13090685885665
83.1615290641785 15.2279701304949
83.5600395202637 23.692954102544
83.9473354816437 23.1025899573632
84.3715252876282 27.8581277108888
84.7385032176971 2.30250379993385
85.1645801067352 1.93256798397498
85.5589830875397 10.5882547311149
85.9375803470612 14.3481842074235
86.3435664176941 9.50740617474453
86.7374198436737 17.0381006099489
87.1536903381348 16.569101178217
87.5287866592407 10.2830296045389
87.927943944931 21.357813519116
88.3459124565125 27.3088723965834
88.736287355423 21.0306229134947
89.1409668922424 21.080097636263
89.5218932628632 34.094267192875
89.9179372787476 29.5693382655294
90.3244364261627 18.6480128347235
90.7179913520813 16.5610712662214
91.1131293773651 9.29881939073529
91.5380389690399 14.7893807055802
91.9342465400696 12.4797340684039
92.3289911746979 22.0565735571605
92.7170422077179 11.6476501982695
93.1170780658722 9.1465865651304
93.5171551704407 24.2263314282982
93.9304356575012 19.6440153851823
94.3490767478943 5.66467854500234
94.7234094142914 1.51319177547885
95.1323461532593 6.05719420390042
95.5392725467682 8.38004150938524
95.9203629493713 10.078675462719
96.3054702281952 16.841484194464
96.707759141922 7.1575900874149
97.1163427829742 8.41087079424795
97.5221474170685 9.62128862770194
97.9280693531036 13.8361959571005
98.3028905391693 15.7786383601612
98.6992218494415 10.4467145041012
99.1074552536011 5.59229158203057
99.4965920448303 6.96212028289865
99.9044208526611 8.61112027916065
100.298995256424 13.2792025817711
100.688905715942 6.95907223475639
101.08878493309 12.2285414699975
101.507484436035 13.3486665235639
};
\end{axis}

\end{tikzpicture}